\documentclass[11pt]{article}

\usepackage{amsfonts,amssymb}

\newcommand{\Z}{\mathbb{Z}}

\newcommand{\F}{\mathbb{F}}

\newcommand{\col}{{\mathcal C}}

\newcommand{\minc}{{\rm mincol}}

\newtheorem{theorem}{Theorem}[section]

\newtheorem{remark}[theorem]{Remark}

\input{epsf.sty}

\date{}

\begin{document}

\title{Minimal Numbers of Fox Colors  and \\ Quandle Cocycle Invariants of Knots}

\author{Masahico Saito 
\\ University of South Florida
}

\maketitle

\begin{abstract}
 Relations will be described between the quandle cocycle 
 invariant and the minimal number of colors used for 
 non-trivial Fox colorings
 of knots and links. In particular, a lower bound for the minimal number is given 
 in terms of the quandle cocycle invariant.
\end{abstract}



\section{Introduction}\label{introsec}

\begin{figure}[htb]
\begin{center}
\mbox{
\epsfxsize=3in
\epsfbox{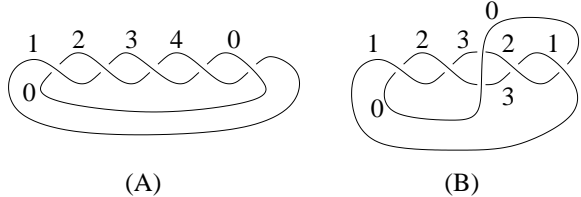}
}
\end{center}
\caption{Teneva's example }
\label{Teneva}
\end{figure}

The determinant of a knot is divisible by a prime $p$ if 
and only if a diagram of the given knot is (Fox) $p$-colorable \cite{Fox}.
Since the knot determinant often takes much larger values 
than the crossing numbers for prime alternating knots with prime determinants,
it was conjectured \cite{HK} that 
for any prime $p$, if an alternating knot $K$ has the determinant $p$,
then any non-trivial coloring of any minimal crossing diagram of $K$ 
assigns distinct colors on its arcs ({\it Kauffman-Harary conjecture}).
The conjecture stays open at the time of writing after extensive 
studies of wide variety of families of knots (see, for example, \cite{APS}).
Considering this situation
the following example is interesting, that 
was discovered by I. Teneva (described  in \cite{HK}):
the $(2,5)$-torus knot $T(2,5)$, which has determinant $5$,
 has a non-alternating, non-minimal diagram 
with only $4$ colors, as depicted in Fig.~\ref{Teneva}(B)
(compare with 
 its minimal alternating diagram  colored by $5$ 
distinct colors as shown in Fig.~\ref{Teneva}(A).)
In \cite{KL}, the minimum number of $p$-colors
$\minc_p(K)$ for knots $K$ was further studied,
and it was proved that if $\minc_p(K)=3$, then $3$ divides the determinant
of a knot $K$. From these facts,
the case $\minc_p(K)=4$ is of interest, and is a focus of this paper.

Quandle cocycle invariants, defined using quandle cohomology theory,
 have been studied   and  
applied  to knots and knotted
surfaces~\cite{CJKLS,CKamS}  in the past several years.
For  Fox $p$-colorings with region colors for knots and links $L$,
quandle cocycle invariant $\Phi_p(L)$  is
defined and  written 
as a multiset (a set with repetitive elements allowed)
in $\Z_p$ (more details will be given in Section~\ref{prelimsec}). 
For example, monochromatic (trivial) 
colorings of a knot contribute copies of $0 \in \Z_p$ 
in $\Phi_p(L)$, but some non-trivial colorings may further contribute $0$ as well.
Hence if $\Phi_p^0 (L) $ denotes the number of zeros in $\Phi_p(L)$,
then  $\Phi_p^0(L)\geq p^2$ for any link $L$,
and 
 $\Phi_p^0(L) > p^2$
if and only if there is a non-trivial $p$-coloring that contributes $0$
to the invariant $\Phi_p(L)$. 
A link is called {\it non-split} if there is no embedded $2$-sphere 
in $3$-space that separates the link non-trivially.
In this paper we prove

\begin{theorem}\label{mainthm}
$(1)$ If  there exists a prime $p>7$ such that 
a non-split link $L$ 
satisfies 
$\Phi_p^0(L)=p^2$, 
 then $\minc_q (L) \geq 5$
for every prime $q>7$.

\noindent
$(2)$ There exist infinitely many  
prime alternating 
  knots $K$ with 
 $\Phi_7^0 (K) =49$,
 and their diagrams that
 are $7$-colored with exactly  $4$ colors.
\end{theorem}

\begin{figure}[htb]
\begin{center}
\mbox{
\epsfxsize=2.5in
\epsfbox{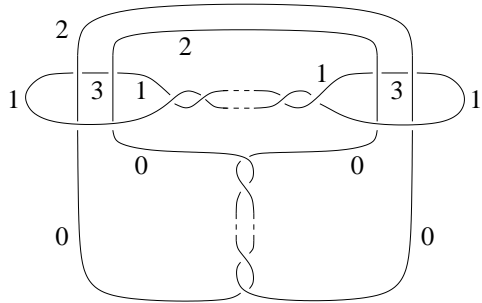}
}
\end{center}
\caption{Links colored with $4$ colors }
\label{linkex}
\end{figure}

\begin{figure}[htb]
\begin{center}
\mbox{
\epsfxsize=1.2in
\epsfbox{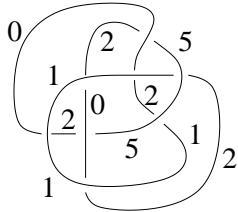}
}
\end{center}
\caption{A $7$ coloring of $5_2$ with $4$ colors}
\label{52}
\end{figure}

For the statement (1),
we exhibit in Fig.~\ref{linkex}  a family of examples of  non-split links 
that are $p$-colorable for any prime $p$,  with exactly  $4$ colors.
The smallest example of the statement (2) of the theorem is 
depicted in Fig.~\ref{52}, which is equivalent to $5_2$, the first 
knot in the table that is $7$-colorable.

In Section~\ref{prelimsec}, we present definitions necessary to prove the main theorem,
and the proof is given in Section~\ref{proofsec}.

\section{Preliminaries}\label{prelimsec}

\begin{figure}[htb]
\begin{center}
\mbox{
\epsfxsize=2in
\epsfbox{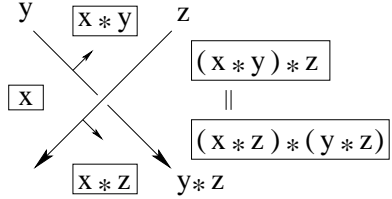}
}
\end{center}
\caption{Coloring rules of arcs and regions}
\label{regionFoxcolors}
\end{figure}

\begin{sloppypar}
Fox $p$-colorings are well known in knot theory,
and their descriptions can be found in \cite{CKamS,K&P,Rolf}, for example.
The definition of quandles, as well as references
  are also  found in \cite{CKamS,K&P}.
In this section we give  brief reviews 
needed
for the statement and proof of Theorem~\ref{mainthm}, and leave details to these references.
The set $\Z_p$ with a binary operation 
$(x, y)\mapsto x*y=2y-x$ mod $p$ for $x, y \in \Z_p$ is called 
the {\it dihedral quandle} of order $p$.
A map from the set of arcs by $\Z_p$ is called a Fox $p$-coloring, 
if it satisfies the condition in Fig.~\ref{regionFoxcolors} for $y, z$ assigned on arcs.
A coloring is regarded as assigning elements of a quandle to arcs,
and the element assigned to an arc is called  the  {\it color} of the arc.
A map from the set of arcs and regions of a planar knot diagram
is called a Fox coloring {\it with region colors},  
if it satisfies the condition in  Fig.~\ref{regionFoxcolors}, where 
the region colors (elements assigned to regions) are framed by squares.
In fact, for defining  Fox colors,  orientations of knots are not necessary.
\end{sloppypar}

Definitions and references for quandle cocycle invariants 
are found in \cite{CKamS}, for example.
Although  orientations are needed to define 
the invariant for general quandles, it is shown by Satoh~\cite{Sat}
that it is well-defined without orientation for dihedral quandles, 
and we use his result for the definition below. 
Let ${\cal C}$ be a Fox $p$-coloring with  region colors of a given diagram
$K$.  
Refer to Fig.~\ref{regionFoxcolors} for the descriptions below.
Select one of the four regions near a  crossing $\tau$, and call it the 
{\it source region}. Let $x_{\tau}$ be the color of the source region.
Select one of the  two under-arcs at $\tau$, 
and call it the {\it source under-arc}.
Let $y_{\tau}$ be the color of  the source under-arc,
and $z_{\tau}$ be the color of the over-arc at $\tau$, 
then $(x_{\tau}, y_{\tau}, z_{\tau})$
is called  the {\it  ordered triple of colors} at a crossing  $\tau$.
Let $n_y$ and $n_z$ be the normal vectors of the source under-arc and the over-arc at $\tau$,
respectively, such that they point from the source region across the arcs to the other regions. 
The local sign of $\tau$ denoted by $\epsilon(\tau)$,
which depends on the choices of the source region 
and the source under-arc, is defined to be positive if the 
ordered vectors $(n_z, n_y)$ agrees with the orientation 
of the plane, and negative otherwise. 
The weight at $\tau$ for the coloring $\col$ is defined by
 $B( {\cal C}, \tau)=\epsilon(\tau) \phi(x_{\tau}, y_{\tau}, z_{\tau})$,
 where $\phi$ is called a {\it quandle $3$-cocycle} described below.
 The quandle ($3$-)cocycle invariant is defined 
 by the multiset $ \Phi_p(K) =\{
\sum_{\tau} B( {\cal C}, \tau ) \ | \ {\cal C}\  \mbox{{\rm ranges over all $p$-colorings
with region colors}}   \} $.
This is a knot invariant, which follows from the properties
satisfied by $\phi$, in particular, the condition called the quandle 
$3$-cocycle condition  (see \cite{CKamS}).
For dihedral quandles, Mochizuki~\cite{Mochi} gave the following explicit formula:
$$\phi(x,y,z)=(x-y)[ (2 z^p - y^p) - (2 z - y )^p ]/p \in \Z_p . $$
It was proved in \cite{Sat} that the weight $B( {\cal C}, \tau)$ does not depend on 
the choices of the source region and the source under-arc when 
Mochizuki's cocycle is used to define the weight. 

As was mentioned in Section~\ref{introsec}, trivial colorings
(those with the same color on all arcs) contribute $0$ to $\Phi_p(K)$.
This is seen from $\phi(x,y,z)$ by setting $y=z$. The contribution is $0$ regardless of region colors, but a color of one region (say, the region at infinity)
uniquely determines colors of all the other regions, 
so that there are at least $p^2$ copies of $0$'s in $\Phi_p(K)$
for any $p$ and $K$. 
Let $\Phi_p^0(K)$ denote the number (multiplicity) of $0$'s in $\Phi_p(K)$,
then we have $\Phi_p^0(K) \geq p^2$, and the inequality is strict 
if and only if there exists a coloring with region colors that is non-trivial 
on arcs, yet contributes $0$. 

For dihedral quandles, we take advantage of the fact that $\Z_p$ forms a field.
Let $p$ be a prime, and $c \in \Z_p$.
If $c \neq 0$, then $c$ is invertible in the field $\Z_p=\F_p$, and
denote by $(1/c)$ the multiplicative inverse of $c$.
For a $p$-coloring of  $\cal$ a knot diagram and $c\in \Z_p$,
 define a coloring $c + \col $  and $c \ \col$ to be  $p$-colorings defined by 
 $ (c + \col )(\alpha)= c + \col (\alpha)$ and 
  $ (c \ \col )(\alpha)= c \cdot \col (\alpha)$, respectively, 
  for every arc $\alpha$ of the diagram.
It is easily seen that they are well-defined.
 Similarly,  if $c \neq 0$, define 
 a coloring $(1/c)\col $ by $((1/c) \col )(\alpha )=(1/c) \cdot \col(\alpha)$
 for every arc $\alpha$ of the diagram.

\section{Proof of Theorem~\ref{mainthm} } \label{proofsec}

Let $L$ be a non-split link.
To prove the statement (1), 
suppose that a diagram $D$ of $L$ has a $p$-coloring $\col$ for a prime $p>7$
with $4$ colors. 
Then we prove that
 the coloring  $\col$ 
can be defined for any prime $q>7$, and makes a trivial contribution ($0 \in \Z_p$)
to the cocycle invariant $\Phi_q(L)$,
 and therefore $\Phi_q(L)>q^2$ for any prime $q>7$.

\bigskip

\begin{figure}[htb]
\begin{center}
\mbox{
\epsfxsize=2.2in
\epsfbox{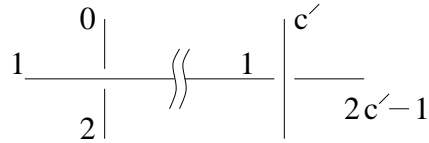}
}
\end{center}
\caption{Specifying $4$ colors}
\label{proof}
\end{figure}

Let $\col(D)$ denote the set of colors that appear in $\col$.
By replacing $\col$ by $(-c) + \col$ for a color $c \in \col(D)$, 
we may assume that $0 \in \col(D)$.
Let $c\neq 0$ be a color that is assigned to  the over-arc 
at the end of an arc colored by $0$, so that at this crossing $\tau$,
one of the under-arcs has color $0$ and the over-arc has color $c$. 
By replacing $\col$ by $(1/c) \col$, we map assume that $c=1$.
Then at $\tau$, the other under-arc is colored by $2$.
The situation is depicted at the left crossing of Fig.~\ref{proof}.

We may assume that there is an arc colored by $1$ that ends at 
an over-arc colored with $c'\neq 1$ at a crossing $\tau'$, since otherwise, 
there are components colored only by $1$ that lie  above all the other components,
and this contradicts the assumption that $L$ is non-split.
Let $c'$ be the color of  the over-arc
of this crossing $\tau'$ that is an end point of an arc colored $1$,
so that the three colors at $\tau'$ 
are $1, c'$, and $2c'-1$.
Since  $\col(D)$ consists of $4$ colors and contains $\{ 0, 1, 2\}$, 
either $c' \in \{ 0, 1, 2\}$ (and we assumed $c'\neq 1$), or $2c'-1  \in \{ 0, 1, 2\}$.
If $c'=0$, then the other under-arc has color $-1$,
and by considering $1+\col$, we may assume that 
$\col(D)=\{0,1,2,3\}$.
If $c'=2$, then the other  under-arc has color $3$,
and again we may assume that $\col(D)=\{0,1,2,3\}$.
If $2c'-1=0$, then $c'=1/2$, and $\col(D)=\{0,1,2,1/2\}$.
If $2c'-1=1$ or $c'$, then $c'=1$, which contradicts our choice of $c' \neq 1$. 
If $2c'-1=2$, then $c'=3/2$, and $\col(D)=\{0,1,2,3/2\}$.

\bigskip 

\noindent 
{\bf Case 1}: $\col(D)=\{0,1,2,3\} \subset \Z_p$.
If there is a crossing at which the over-arc is colored by $0$, 
then the under-arcs are colored by unordered 
pairs of elements of $\Z_p$ $(0,0)$, $(1,-1)$, $(2,-2)$ or $(3,-3)$,
and only $(0,0)$ is possible since $p>5$. 
For an over-arc colored $3$, 
 possible colors of the under-arcs are 
$(0,6)$, $(1,5)$, $(2,4)$ and $(3,3)$,  all numbers considered in $\Z_p$, and 
only $(3,3)$ is possible since $p>5$.
Hence  $0$ and $3$ cannot be a color of an over-arc other than at 
a trivially colored crossing.

If $1$ is at an over-arc, then the possible colors of under-arcs are 
$(0,2)$, $(1,1)$ and  $(3, -1)$, the last of which is impossible.
For $2$ at an over-arc, then  possible colors of the under-arcs are
$(0,4)$, $(1,3)$ and $(2,2)$, the first of which is impossible.
In summary, 
the possible colors at a crossing are, other than a constant coloring at a crossing, 
the over-arc $1$, under-arcs $(0,2)$ or 
the over-arc $2$, under-arcs $(1,3)$. 
Note that  the arc colored $0$ always ends at an over-arc colored $1$, and 
$3$ ends at $2$.
Note also that this coloring rule by $\{ 0,1,2,3\}$ is valid for any $p$.

\bigskip

At this point we know that 
there is no knot that is colored by this pattern, since 
the diagram is non-trivially colored 
for any $p$, and the existence of such 
colorings  would imply~\cite{Fox}  that the determinant of the knot is divisible 
by all primes $p$, a contradiction (knot determinant takes values in odd integers,
see \cite{Rolf}, for example).

\bigskip

\begin{figure}[htb]
\begin{center}
\mbox{
\epsfxsize=1.5in
\epsfbox{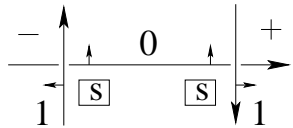}
}
\end{center}
\caption{Pair crossings' contributions cancel}
\label{pair}
\end{figure}

Now we evaluate the contribution to the cocycle invariant  of this coloring
with an arbitrary fixed region colors for the link $L$.
Let $(D_0, \col_0)$ be this  coloring  by $\{ 0,1,2,3\}$ of the given diagram, 
and let $B(D_0, \col_0)$ be the evaluation by the Mochizuki $3$-cocycle
(which is the contribution to the cocycle invariant $\Phi_p (L)$ of 
this coloring).
Let $D_1$ be the diagram obtained from $D_0$ by performing  
(arbitrary) smoothings
at all crossings colored by constant colors $0, 1, 2$ or $3$, and give 
 the inherited coloring $\col_1$
to obtain $(D_1, \col_1)$.
 The link type of $D_1$ may no longer be $L$, but we retain 
 $B(D_1, \col_1)=B(D_0, \col_0)$.
 Every arc colored $0$ or $3$ of the diagram $D_1$ has
 no crossing with constant (trivial) colors after smoothings,
 and they end  at over-arcs with distinct colors.
 Every crossing has a single under arc colored by either $0$ or $3$,
 exclusively. Hence the set of arcs colored by $0$ or $3$ gives
 rise to  pairings of crossings, by declaring that 
 the end points of an arc colored by $0$ or $3$ 
 are paired.  The contribution to the  
 weight $B(D_1, \col_1)$ of each pair is zero. 
 This is seen by selecting a region shared by the pair of crossings as 
 the source region, and the shared under-arc as the source under-arc.
 Then  the crossings have the same 
  ordered triple of colors and opposite local signs, 
  see Fig.~\ref{pair}.
 Hence we have  $B(D_1, \col_1)=B(D_0, \col_0)=0$, a contradiction
 to the assumption that $\Phi_p(L)=p^2$.
 The argument is valid for all prime $q>7$.

\noindent 
{\bf Case 2}: $\col(D)=\{0,1,2,1/2\} \subset \Z_p$. We may assume that these are distinct, 
and that $1/2 \neq 3$.
If there is a crossing at which the over-arc is colored by $0$, 
then possibilities of colors of the under-arcs are unordered
pairs $(0,0)$, $(1, -1)$, $(2, -2)$ and $(1/2,-1/2)$.
All colors must be from $\col(D)=\{0,1,2,1/2\}$.
Since  $p>5$ we have $\{ -1,  -2\} \neq \{ 1, 2 \}$, as well as
 $1/2 \neq -1/2$,
$-1/2 \neq 1$,  $ -1\neq 1/2$,
$ -2 \neq 1/2$.
Hence only $(0,0)$ is possible.
For an over-arc colored $1$, the under-arcs are colored
$(0,2)$, $(1,1)$, or $(1/2,3/2)$, but $3/2$ cannot be any of
$\col(D)=\{0,1,2,1/2\}$, so only $(0,2)$ and $(1,1)$ are possible.
If $2$ is at an over-arc, then the under-arcs may be  colored 
$(0,4)$, $(1,3)$, $(2, 2)$ or $(3/2,5/2)$, 
and again only $(2,2)$ is possible if $p>7$, 
and $(0,4)$ is possible for $p=7$ where $\col(D)=\{0,1,2,1/2=4\}$.
For $1/2$ at an over-arc, the under-arcs are colored 
$(0,1)$ or $(1/2,1/2)$.
In summary, 
the possible colors at a crossing are, other than the constant coloring at a crossing, 
the over-arc $1$, under-arcs $(0,2)$,
and 
the over-arc $1/2$, under-arcs $(0,1)$. 
Since $1/2$ does not appear on an under-arc other than 
crossings colored by a constant $1/2$, the components colored by $1/2$
lies above all the others, and contradicts that $L$ is non-split. 
Hence this case does not occur for $p>7$.

\bigskip
 
\noindent 
{\bf Case 3}: $\col(D)=\{0,1,2,3/2\} \subset \Z_p$. We may assume that these are distinct, 
and that $3/2 \neq 3$.
If there is a crossing at which the over-arc is colored by $0$, 
then possible colors of the under-arcs are
pairs $(0,0)$, $(1,-1)$, $(2,-2)$ and $(3/2,-3/2)$.
Since  $p>5$ we have $\{-1, -2\} \neq \{ 1, 2 \}$, $3/2 \neq -3/2$,
$-3/2 \neq 0, 1$, and $p-1\neq 3/2$.
Hence $(1,-1)$ and $(3/2,-3/2)$ are impossible.
If $-2=3/2$, or $-3/2=2$, then $p=7$.
Hence $(1, -1)$ is impossible for any $p>5$ and 
$(2,-2)$  and $(3/2,-3/2)$  are possible only for $p=7$, in which case both are $(5, 2)$.
For over-arc colored $1$, the under-arcs are colored
$(0,2)$, $(1,1)$ or $(3/2,1/2)$, 
and $(3/2,1/2)$ is impossible as 
$1/2\neq 0, 1, 2, 3/2$ 
for odd prime $p>5$.
If $2$ is at an over-arc, then the under-arcs are colored 
$(0,4)$, $(1,3)$, $(2, 2)$ or $(3/2,5/2)$, 
and only $(2,2)$ is possible since $p>5$.
For $3/2$ at an over-arc, the under-arcs are colored 
$(0,3)$, $(1,2)$ or $(3/2,3/2)$,
and $(0,3)$ is impossible.
In summary, 
the possible colors at a crossing are, other than a constant coloring at a crossing, 
the over-arc $0$, under-arcs $(2,5)$ if $p=7$,
the over-arc $1$, under-arcs $(0,2)$ and 
the over-arc $3/2$, under-arcs $(1,2)$. 
Note that  $2$ cannot be a color of an over-arc other than at 
a trivially colored crossing.
If $p>7$, then  $3/2$ does not appear as a color of an under-arc,
and therefore the components colored by $3/2$ splits,
and it is a contradiction. 
Hence the proof is complete for $p>7$.

\begin{figure}[htb]
\begin{center}
\mbox{
\epsfxsize=2in
\epsfbox{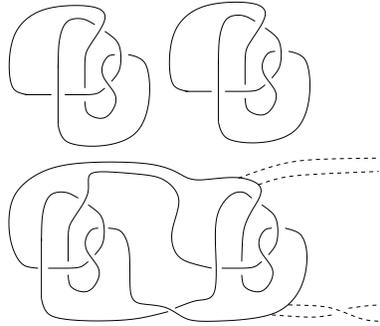}
}
\end{center}
\caption{Connecting  $5_2$}
\label{connect}
\end{figure}

For $p=7$, an example of a knot that satisfies the condition 
stated in the theorem is already exhibited in Fig.~\ref{52}. 
It remains to construct an infinite family of such.
In Fig.~\ref{connect}, a construction is illustrated to obtain such an 
infinite family by connecting any finite number of copies of $5_2$. 
{}From the construction, it is reduced alternating, hence they are all
distinct for distinct numbers of copies used, from \cite{Mura,Th}.
The connection is made at arcs with the same color in Fig.~\ref{52}, 
and by going back to copies of diagrams in Fig.~\ref{52}, the colorings
with only four colors are obtained.
Finally, we observe that for infinitely many of them, 
 every non-trivial coloring contributes 
non-trivial values to the cocycle invariant.
Let $K_n$ be the knot constructed above with $n$ copies of $5_2$.
Any non-trivial coloring of $5_2$ contributes 
$1$, $2$, or $4 \in \Z_7$ 
to the invariant $\Phi_7(K_n)$.
Then for each  contribution, 
the contribution of the induced coloring for $K_n$ is 
$n$, $2n$, and $4n \in \Z_7$, 
respectively, so that the 
invariant is non-trivial for all $n$ not divisible by $7$.
This provides an infinite family as required.
$\Box$

\begin{remark}{\rm
It was proved in \cite{KO}, in fact, that any $7$-colorabe knot has a diagram with exactly four colors.
Using her result,  Theorem~\ref{mainthm} (2) can be proved simply by providing 
infinitely many $7$-colorable knots $K$ with $\Phi_p^0 (K)=49$.
} \end{remark}

\begin{remark}{\rm
In terms of quandle homology theory \cite{CKamS01},
Theorem\ref{mainthm} can be restated as follows:
(1) For any $p>7$, any $p$-coloring of a link diagram with $4$ colors represents 
a null-homologous $3$-cycle in $H_3^{\rm Q}(\Z_p, \Z_p)$.
(2) There 
 exist infinitely many  
prime alternating 
  knots $K$, and $7$-colorings 
 of their diagrams with exactly  $4$ colors,
 that represent non-zero homology classes 
 of $H_3^{\rm Q}(\Z_7, \Z_7)$. 
The minimal number of elements of a quandle used to represent 
cycles that are non-trivial in quandle homology is of interest
from this view point.
}
\end{remark}

\section*{Acknowledgments}
The author was supported in part by NSF Grant DMS 
 \#0603876.
He is thankful to J.\,Scott Carter, L.H.\,Kauffman, P.\,Lopes and K.\,Oshiro for 
 valuable conversations and comments.

\end{document}